\documentclass{article}
\usepackage{amssymb}
\usepackage{mathrsfs}
\usepackage{cite}

\usepackage{stmaryrd}
\usepackage{amsfonts}
\usepackage{amsmath}
\usepackage{bm}
\usepackage{indentfirst}
\usepackage{array}
\usepackage{arydshln}
\usepackage{graphicx}
\usepackage{listings}
\usepackage{epstopdf}
\numberwithin{equation}{section}
\usepackage{geometry}
\usepackage{amsthm}
\newtheorem{theorem}{Theorem}

\newtheorem{lemma}{Lemma}
\newtheorem{rem}{Remark}
\newtheorem{proposition}{Proposition}
\title{\Large \bf{Decay properties of solutions toward shock waves of the scalar conservation law with linear and $p$-Laplacian viscosity}
\setcounter{footnote}{-1}
\author{Yechi Liu$^\dag$
\noindent\footnote{\dag\quad E-mail: lyc9009@sina.cn.}\\
\small College of Science, National University of Defence Technology, Changsha 410003, P.R.China}}
\date{}
\begin{document}\large
\maketitle
\textbf{Abstract}. In this paper, we discuss the asymptotic behaviour of weak solutions to the Cauchy problem toward the viscous shock waves for the scalar viscous conservation law. We firstly consider the case that the flux function is the quadratic Burgers flux and obtain the time-decay rate in $L^\infty$-norm for the cases with degenerate $p$-Laplacian viscosity and with linear viscosity, respectively. Moreover, we also give the time-decay rate in $L^\infty$-norm with general flux and linear viscosity. All of these results do not involve smallness conditions on the initial data.

\vskip 0.2in


\section{\Large Introduction and main results}
Consider the one-dimensional scalar conservation law with a nonlinear degenerate viscosity
\begin{equation}\label{0}
u_t+f(u)_x=\mu\big(|u_x|^{p-1}u_x\big)_x,\qquad t>0,x\in\mathbb R
\end{equation}
with initial data
\begin{equation}\label{inid}
u(0,x)=u_0(x),\qquad x\in\mathbb R
\end{equation}
and
\begin{equation}\label{bndd}
\lim_{x\rightarrow\pm\infty}u(t,x)=u_\pm,\qquad t\geqslant0.
\end{equation}
Here, the flux $f\in C^2(\mathbb R)$ satisfies $f(0)=f^\prime(0)=0$ and $f^{\prime\prime}(v)\geqslant c_f>0$ for any $v\in\mathbb R$. The constants $u_\pm$ are prescribed far field states and $\mu>0$. By a scaling
\begin{equation*}
t\rightarrow\mu\tau,\quad x\rightarrow\mu^\frac{2}{p+1}y
\end{equation*}
we can obtain that
\begin{equation*}
u_\tau+\mu^\frac{p-1}{p+1}f(u)_y=\big(|u_y|^{p-1}u_y\big)_y.
\end{equation*}
Thus, without loss of generality, we can suppose $\mu=1$ in the rest of this paper. It is known that if $p=1$ and $f(u)=\frac{1}{2}u^2$, \eqref{0} becomes the viscous Burgers equation
\begin{equation}\label{be}
u_t+uu_x=u_{xx}.
\end{equation}
The linear viscosity term $u_{xx}$ stands for Newtonian fluid. The nonlinear viscosity term $\big(|u_x|^{p-1}u_x\big)_x$ stands for non-Newtionian fluid, such as blood, honey, butter, whipped cream, etc. (see \cite{yos15}). Such a viscosity term is also called Ostwald-de Waele type viscosity. For more details, see \cite{yos17,yos17a}.

We also note that our model \eqref{0} is motivated by the following Lady\v{z}enskaja model
\begin{equation*}
\partial_tu_i+\sum_ju_j\partial_ju_i+\partial_iP=\sum_j\partial_j\big((\mu_0+\mu_1|\nabla u|^{p-1})\partial_ju_i\big)+f_i,
\end{equation*}
where $|\nabla u|=\big(\sum_{i,j}|\partial_iu_j|^2\big)^\frac{1}{2}$. This is the incompressible Navier-Stokes equation with the power-law type nonlinear viscosity (see \cite{du91,lad70,mal06}). For one-dimensional case, the viscosity becomes
\begin{equation*}
\mu_0u_{xx}+\mu_1\big(|u_x|^{p-1}u_x\big)_x.
\end{equation*}
When $\mu_0=0$, it is the viscosity in \eqref{0}.

We are interested in the asymptotic behavior and precise estimates in time of the global solution to the problem (\ref{0}-\ref{bndd}). It can be expected that the large time behavior is closely related to the solution of the corresponding Cauchy problem
\begin{equation}\label{rp}
u_t+f(u)_x=0
\end{equation}
with Riemann initial data
\begin{equation}\label{rid}
u_0(x)=\left\{\begin{aligned}u_-,\quad x<0,\\u_+,\quad x>0,\end{aligned}\right.
\end{equation}
which is also known as Riemann problem. Asymptotic results have a long history starting with the paper of Il'in-Ole\u{\i}nik \cite{ili60,ili64} for the conservation law with $p=1$, i.e.,
\begin{equation}\label{lvp}
u_t+f(u)_x=u_{xx},
\end{equation}
and $f$ being genuinely nonlinear (i.e. $f^{\prime\prime}\neq0$). They showed that if $f^{\prime\prime}(u)>0(u\in\mathbb R)$ and $u_-<u_+$ (or $f^{\prime\prime}(u)<0(u\in\mathbb R)$ and $u_->u_+$), which implies the Riemann solution of \eqref{rp} and \eqref{rid} consists of a single rarefaction wave solution, the global in time solution of \eqref{lvp}, \eqref{inid} and \eqref{bndd} tends toward the rarefaction wave as $t\rightarrow\infty$; if $f^{\prime\prime}(u)>0(u\in\mathbb R)$ and $u_->u_+$ (or $f^{\prime\prime}(u)<0(u\in\mathbb R)$ and $u_-<u_+$), which implies the Riemann solution of \eqref{rp} and \eqref{rid} consists of a single shock wave, the global in time solution of \eqref{lvp}, \eqref{inid} and \eqref{bndd} approaches to a corresponding smooth traveling solution, which is called viscous shock wave, with a spacial shift. Hattori-Nishihara \cite{hat91} also gave the asymptotic time-decay rate $(1+t)^{-\frac{1}{2}(1-\frac{1}{q})}$ in $L^q$-norm ($1\leqslant q\leqslant\infty$) for the solution toward the single rarefaction wave. For viscous shock waves, Kawashima-Matsumura \cite{kaw85} gave a decay result, which says that algebraic decay at $\pm\infty$ in space is transformed to algebraic decay of the perturbation in time. On the other hand, Sattinger \cite{sat76} obtained a result similar to the decay theorem of Il'in-Ole\u{\i}nik, but without the assumption of convexity on the flux. Jones-Gardner-Kapitula \cite{jon93} also obtained the stability of shock waves for the case with any $C^2$ nonlinear flux, provided that the system supports a viscous profile, by using the semigroup method. Recently, Huang-Xu \cite{hua22} got the decay rate to the viscous shock wave in $L^\infty$-norm.

Note that all above results on shock waves are based on the smallness of initial data. For large initial data, Freist\"{u}hler-Serre \cite{fre98} proved the $L^1$ stability of viscous shock waves by combining energy estimates, a lap-number argument and a specific geometric observation on attractor of steady states. Then, Nishihara-Zhao \cite{nis02} further obtained the convergence rate toward the viscous shock waves in $L^\infty$-norm by combining the idea of \cite{kaw85} with the result of \cite{fre98} under some restriction conditions on initial data. In addition, Kang-Vasseur \cite{kan17} showed the decay rate in $L^2$-norm in the context of small perturbations of the quadratic Burgers flux. In this paper, we use the area inequality to obtain the time-decay rate without any assumption on the spatial decay of the initial data.

When the flux function $f$ is not uniformly genuinely nonlinear, there are also some results, see \cite{mat94,mat12,yos14} for example.

For the case with $p>1$, there are few results on the asymptotic behavior for the problem (\ref{0}-\ref{bndd}). When the flux function $f$ is genuinely nonlinear on $\mathbb R$ and $f^{\prime\prime}>0$, Matsumura-Nishihara \cite{mat94a} showed that if the far field states are same (i.e. $u_-=u_+$), the solution to (\ref{0}-\ref{bndd}) tends toward the constant state $u_\pm$, and if $u_-<u_+$, the solution tends toward a single rarefaction wave. Yoshida \cite{yos15} gave the precise time-decay rate in $L^q$-norm ($1\leqslant q\leqslant\infty$) for the asymptotics of the results in \cite{mat94a}, depending on the space of the initial function $u_0$. When the flux function $f$ is genuinely nonlinear on $\mathbb R$ except a finite interval $I$ and linearly degenerate on $I$ (i.e. $f^{\prime\prime}=0$ on $I$), Yoshida \cite{yos17} obtained the asymptotics and its time-decay estimates of the solution to (\ref{0}-\ref{bndd}) with $u_-<u_+$ tends toward the multiwave pattern of the combination of the viscous contact wave and the rarefaction wave. For viscous shock wave ($u_->u_+$), however, there is no result yet. In this paper, we are motivated by Kang-Vasseur \cite{kan17} and obtain the time-decay rates in $L^2$-norm and $L^\infty$-norm to viscous shock waves with large initial data.

\vskip 0.2in

A bounded and measurable function $u(t,x)$ with $u_x\in L_{\rm loc}^{p}\big([0,\infty)\times\mathbb R\big)$ is called a weak solution to (\ref{0}-\ref{bndd}), if
\begin{equation*}
\int_0^T\!\!\int_\mathbb R\Big(u\zeta_t+\big(f(u)-|u_x|^{p-1}u_x\big)\zeta_x\Big){\rm d}x{\rm d}t+\int_\mathbb Ru_0(x)\zeta(0,x){\rm d}x=0
\end{equation*}
holds for any $T>0,\zeta\in C_0^1\big([0,T)\times\mathbb R\big)$ and
\begin{equation*}
\lim_{x\rightarrow\pm\infty}u(t,x)=u_\pm.
\end{equation*}

In this paper, we suppose $u_->u_+$ and firstly consider the case, in which the flux function is the Burgers flux and the viscosity is the degenerate $p$-Laplacian type with $p>1$, that is
\begin{equation}\label{1}
\left\{\begin{aligned}
&u_t+uu_x=\big(|u_x|^{p-1}u_x\big)_x,\\
&u(0,x)=u_0(x),\\
&\lim_{x\rightarrow\pm\infty}u(t,x)=u_\pm.
\end{aligned}\right.
\end{equation}
There is a travelling wave solution $U(\xi)=U(x-\gamma t)$ of \eqref{0} with $\xi=x-\gamma t$, connecting $u_-$ at $x=-\infty$ to $u_+$ at $x=+\infty$, and satisfying $U^\prime\leqslant0$ and
\begin{equation}\label{U}
\left\{\begin{aligned}
&\big((-U^\prime)^p+\frac{1}{2}U^2-\gamma U\big)^\prime=0,\\
&\lim_{\xi\rightarrow\pm\infty}U(\xi)=u_\pm, \quad \lim_{\xi\rightarrow\pm\infty}U^\prime(\xi)=0,
\end{aligned}\right.
\end{equation}
where $\gamma=\frac{u_-+u_+}{2}$ is the speed of the viscous shock wave determined by the Rankine-Hugoniot Condition
\begin{equation*}
\gamma(u_--u_+)=\frac{1}{2}u_-^2-\frac{1}{2}u_+^2.
\end{equation*}
Some properties of $U$ will be given in Section 2.

We then have the following theorem.
\begin{theorem}
If $p>1$, $u_0-U\in L^2(\mathbb R)$ and $u_{0x}\in L^{p+1}(\mathbb R)$, then \eqref{1} has a unique weak solution $u(t,x)$ satisfying
\begin{gather*}
u-U\in C\big([0,\infty);L^2(\mathbb R)\big)\cap L^\infty\big([0,\infty);L^2(\mathbb R)\big),\\
u_x\in L^\infty\big([0,\infty);L^{p+1}(\mathbb R)\big)\cap L^{p+1}\big([0,\infty)\times\mathbb R\big),\\
\big(|u_x|^{p-1}u_x\big)_x\in L^2\big([0,\infty)\times\mathbb R\big).
\end{gather*}
Furthermore, there is a $p_0\in(\frac{39}{20},\frac{59}{30})$ such that, if $1<p\leqslant p_0$ and $u_0-U\in L^1(\mathbb R)$, it holds
\begin{equation}\label{re0}
\big\|u\big(t,\cdot+X(t)\big)-U(\cdot)\big\|_{L^2(\mathbb R)}\leqslant C(1+t)^{-\frac{1}{4p}},
\end{equation}
and hence,
\begin{equation*}
\big\|u\big(t,\cdot+X(t)\big)-U(\cdot)\big\|_{L^\infty(\mathbb R)}\leqslant C(1+t)^{-\frac{1}{2p(p+3)}},
\end{equation*}
where the shift $X(t)$ satisfies
\begin{equation}\label{X}
\left\{\begin{aligned}
&X^\prime(t)=\gamma-\frac{1}{2(u_--u_+)}\int_\mathbb R\Big(u\big(t,x+X(t)\big)-U(x)\Big)U^\prime(x){\rm d}x,\\
&X(0)=0
\end{aligned}\right.
\end{equation}
and $p_0$ will be given by the proof of Lemma \ref{abm} in Appendix B.
\end{theorem}

\vskip 0.2in

Next, if $p=1$ in \eqref{1}, which means
\begin{equation}\label{2}
\left\{\begin{aligned}
&u_t+uu_x=u_{xx},\\
&u(0,x)=u_0(x),
\end{aligned}\right.
\end{equation}
we do not need to suppose the boundary condition \eqref{bndd} holds. The existence of solution to \eqref{2} has been proved, and we obtain the following result on the decay rate.
\begin{theorem}
If $u_0-U\in L^1(\mathbb R)\cap L^\infty(\mathbb R)$, then the solution to \eqref{2} satisfies
\begin{gather*}
\|u\big(t,\cdot+X(t)\big)-U(\cdot)\|_{L^2(\mathbb R)}\leqslant C(1+t)^{-\frac{1}{4}},\\
\|u\big(t,\cdot+X(t)\big)-U(\cdot)\|_{L^\infty(\mathbb R)}\leqslant C(1+t)^{-\frac{1}{6}},
\end{gather*}
where the viscous shock wave $U$ is the solution to \eqref{U} with $p=1$ and the shift $X(t)$ is given by \eqref{X}.
\end{theorem}
\begin{rem}
In Theorem 2, since we do not have the condition of $u_{0x}\in L^{p+1}(\mathbb R)$, it fails to obtain that $u_x\in L^{p+1}(\mathbb R)$, so we will use the linearity of the viscosity to obtain the $L^2$ decay rate of the first-order derivative. This method, however, is not applicable in the case with degenerate viscosity.
\end{rem}

In addition, we consider the general flux function, which means
\begin{equation}\label{3}
\left\{\begin{aligned}
&u_t+f(u)_x=u_{xx},\\
&u(0,x)=u_0(x).
\end{aligned}\right.
\end{equation}
Since the method given in \cite{kan17} is not applicable for general convex flux, we obtain the convergence result in another way, that is
\begin{theorem}
Assume $f\in C^2(\mathbb R)$, $f(0)=f^\prime(0)=0$ and $f^{\prime\prime}(v)\geqslant c_f>0$ for any $v\in\mathbb R$. If $u_0-U\in L^1(\mathbb R)\cap L^\infty(\mathbb R)$ and $\int_{-\infty}^x(u_0-U){\rm d}y\in L^2(\mathbb R)$, then there exists a space shift $y$ such that the solution to \eqref{3} satisfies
\begin{equation}\label{bu10}
\begin{aligned}
\|u\big(t,\cdot\big)-\tilde U(\cdot+y+\gamma t)\|_{L^2(\mathbb R)}\leqslant C(\delta)(1+t)^{-\frac{1}{8}+\delta},\\
\|u\big(t,\cdot\big)-\tilde U(\cdot+y+\gamma t)\|_{L^\infty(\mathbb R)}\leqslant C(\delta)(1+t)^{-\frac{1}{6}+\delta},
\end{aligned}
\end{equation}
where $\delta>0$ is any small constant, $\tilde U$ is the viscous shock wave with the general flux $f$.
\end{theorem}
\begin{rem}
For the existence and properties of the viscous shock wave $\tilde U$, see \cite{mat94}.
\end{rem}

\vskip 0.2in

The rest of this paper is organized as follows. In Section 2, we give some properties of $U$ and some lemmas. Section 3 is devoted to the proof of time-decay rates in Theorems 1-3. For Theorems 1 and 2, a special weighted Poincar\'{e}'s inequality is used and the proof of Theorem 3 is based on the $L^1$ result in \cite{fre98} and energy estimate progress in \cite{hua22}. The proof of existence part in Theorem 1 and some lemmas are given in Appendix.

\vskip 0.2in

\noindent\textbf{Notations}. For function spaces, $L^q=L^q(\mathbb R)$ denotes the usual Lebesgue space with norm $\|\cdot\|_q$, which means
\begin{equation*}
\|v\|_q=:\left(\int_\mathbb R|v(x)|^q\,\textrm{d}x\right)^\frac{1}{q},\quad 1\leqslant q<\infty
\end{equation*}
and
\begin{equation*}
\|v\|_\infty=:{\rm esssup}|v|.
\end{equation*}
In addition, we use $c$ and $C$ to represent uncertain positive constants suitably small and large respectively. In particular, $c(a_1,a_2,\cdots)$ and $C(b_1,b_2,\cdots)$ represent that the constant $c$ and $C$ depend only on $a_1,a_2,\cdots$ and $b_1,b_2,\cdots$, respectively.

\section{\Large Prelimilaries}
Firstly, we will give some properties of the viscous shock wave $U$. Let
\begin{equation}\label{cu}
C_U=(-U^\prime)^p+\frac{1}{2}U^2-\gamma U,
\end{equation}
then, from \eqref{U}, we have
\begin{equation*}
C_U-\frac{1}{2}u_\pm^2+\gamma u_\pm=\lim_{\xi\rightarrow\pm\infty}\big(U^\prime(\xi)\big)^p=0,
\end{equation*}
which means
\begin{equation*}
C_U=\frac{1}{2}u_\pm^2-\gamma u_\pm.
\end{equation*}
It then follows from \eqref{cu} that
\begin{equation}\label{upr}
U^\prime=-\big(\frac{1}{2}(u_--U)(U-u_+)\big)^\frac{1}{p}.
\end{equation}
Formally, for any $\xi\in\mathbb R$, from \eqref{upr}, it holds
\begin{equation*}
\xi=-\int_{U(0)}^{U(\xi)}\frac{{\rm d}w}{\big(\frac{1}{2}(u_--U)(U-u_+)\big)^\frac{1}{p}}.
\end{equation*}
Without loss of generality, we suppose $U(0)=\frac{u_-+u_+}{2}=:U_0$.

Noting that $w\in(u_+,u_-)$, set
\begin{equation*}
x_L=-\int_{U_0}^{u_-}\frac{{\rm d}w}{\big(\frac{1}{2}(u_--U)(U-u_+)\big)^\frac{1}{p}},\quad x_R=-\int_{U_0}^{u_+}\frac{{\rm d}w}{\big(\frac{1}{2}(u_--U)(U-u_+)\big)^\frac{1}{p}},
\end{equation*}
then $-\infty<x_L<0<x_R<+\infty$ from $p>1$ and $x_L=-\infty,x_R=+\infty$ for $p=1$. Therefore, $U(\xi)$ satisfies
\begin{equation*}
\xi=-\int_{U_0}^{U(\xi)}\frac{{\rm d}w}{\big(\frac{1}{2}(u_--U)(U-u_+)\big)^\frac{1}{p}}
\end{equation*}
for $\xi\in(x_L,x_R)$ \big(hence $u_+<U<u_-$, $U^\prime(\xi)<0$ for $x_L<\xi<x_R$\big) and
\begin{equation*}
U(\xi)=\left\{\begin{aligned}
&u_-,\quad \xi\leqslant x_L,\\
&u_+,\quad \xi\geqslant x_R.
\end{aligned}\right.
\end{equation*}
Furthermore, from \eqref{upr}, it holds
\begin{equation}\label{uprest}
|U^\prime(\xi)|\leqslant2^{-\frac{3}{p}}(u_--u_+)^\frac{2}{p}
\end{equation}
for any $\xi\in\mathbb R$.

\vskip 0.2in

Next, we need the existence and uniqueness of the solution to \eqref{1}, that is
\begin{proposition}\label{ext}
If $p>1$, $u_0-U\in L^2(\mathbb R)$ and $u_{0x}\in L^{p+1}(\mathbb R)$, then \eqref{1} has a unique weak solution $u(t,x)$ satisfying
\begin{gather*}
u-U\in C\big([0,\infty);L^2(\mathbb R)\big)\cap L^\infty\big([0,\infty);L^2(\mathbb R)\big),\\
u_x\in L^\infty\big([0,\infty);L^{p+1}(\mathbb R)\big)\cap L^{p+1}\big([0,\infty)\times\mathbb R\big),\\
\big(|u_x|^{p-1}u_x\big)_x\in L^2\big([0,\infty)\times\mathbb R\big).
\end{gather*}
\end{proposition}
\noindent The proof of Proposition \ref{ext} will be given in Appendix A.

\vskip 0.2in

To obtain the time-decay rate, we also need the following Lemmas.
\begin{lemma}\label{yt}
Let $y(t)$ be a non-negative differentiable function on $[0,\infty)$ and $y^\prime(t)\in L^1_{\rm loc}\big([0,\infty)\big)$, such that
\begin{equation}\label{yt0}
y^\prime(t)+a(1+t)^\alpha y(t)^{1+\beta}\leqslant b(1+t)^{-\gamma}
\end{equation}
for $t\geqslant s$ with some $s\geqslant0$, where $a,\beta,\gamma>0,\alpha,b\geqslant0$ are constants independent of $t$. Then,
\begin{equation}\label{yt1}
y(t)\leqslant C_0(1+t)^{-\mu}
\end{equation}
for any $t\geqslant s$, in which
\begin{equation*}
\mu=\min\left\{\frac{\alpha+\gamma}{1+\beta},\frac{1+\alpha}{\beta}\right\},\quad C_0=\max\left\{\left(\frac{2b}{a}\right)^\frac{1}{1+\beta},
\left(\frac{2(1+\alpha)}{a\beta}\right)^\frac{1}{\beta}\right\}.
\end{equation*}
\end{lemma}
\begin{proof}
For simplicity, we choose $s=0$ in the proof. We firstly consider the following inequality on $t\in[t_0,\infty)$
\begin{equation}\label{yt2}
z^\prime(t)+\tilde a(1+t)^\alpha z(t)^{1+\beta}\leqslant0
\end{equation}
with any $t_0\geqslant0$, in which the constant $\tilde a>0$. Solving \eqref{yt2}, we have
\begin{equation*}
z(t)\leqslant\left(z(t_0)^{-\beta}
+\frac{\tilde a\beta}{1+\alpha}\big((1+t)^{1+\alpha}-(1+t_0)^{1+\alpha}\big)\right)^{-\frac{1}{\beta}}.
\end{equation*}
If
\begin{equation*}
z(t_0)^{-\beta}-\frac{\tilde a\beta}{1+\alpha}(1+t_0)^{1+\alpha}>0,
\end{equation*}
then
\begin{equation*}
z(t)\leqslant\left(\frac{1+\alpha}{\tilde a\beta}\right)^\frac{1}{\beta}(1+t)^{-\frac{1+\alpha}{\beta}}
\end{equation*}
for any $t\geqslant t_0$. Otherwise,
\begin{equation*}
z(t_0)\leqslant\left(\frac{1+\alpha}{\tilde a\beta}\right)^\frac{1}{\beta}(1+t_0)^{-\frac{1+\alpha}{\beta}}.
\end{equation*}

Next, suppose \eqref{yt1} fails, which means there exists $t_1>0$ such that
\begin{equation}\label{yt3}
y(t)>\left(\frac{2b}{a}\right)^\frac{1}{1+\beta}(1+t)^{-\frac{\alpha+\gamma}{1+\beta}}
\end{equation}
and
\begin{equation}\label{yt4}
y(t)>\left(\frac{2(1+\alpha)}{a\beta}\right)^\frac{1}{\beta}(1+t)^{-\frac{1+\alpha}{\beta}}
\end{equation}
for $t_1\leqslant t\leqslant t_1+\delta$ with some small $\delta>0$. From \eqref{yt3} we have
\begin{equation*}
b(1+t)^{-\gamma}<\frac{a}{2}(1+t)^\alpha y(t)^{1+\beta},
\end{equation*}
which, together with \eqref{yt0} implies
\begin{equation*}
y(t)^\prime+\frac{a}{2}(1+t)^\alpha y(t)^{1+\beta}<0.
\end{equation*}
Applying the discussion about \eqref{yt2} above, we obtain
\begin{equation*}
y(t)\leqslant\left(\frac{2(1+\alpha)}{a\beta}\right)^\frac{1}{\beta}(1+t)^{-\frac{1+\alpha}{\beta}},
\end{equation*}
which contradicts \eqref{yt4}. Thus, \eqref{yt1} holds true.
\end{proof}
\begin{lemma}\label{abm}
There is a $p_0\in(\frac{39}{20},\frac{59}{30})$ such that for $1\leqslant p<p_0$,
\begin{equation}\label{ab}
\big(|a|^{p-1}a-|b|^{p-1}b\big)(a-b)\geqslant\frac{5}{6}\Big(c(p)|a-b|^{p-1}
+\big(\!\max\{|a|,|b|\}\big)^{p-1}\Big)(a-b)^2
\end{equation}
for any $a,b\in\mathbb R$, where $c(p)$ is a constant only depending on $p$.
\end{lemma}
Lemma \ref{abm} will be proved in Appendix B.
\begin{lemma}[\!\cite{fre98}]\label{phil1est}
If $u_0-U\in L^1(\mathbb R)$ and $\int_\mathbb R\big(u_0(x)-U(x)\big){\rm d}x=0$, then the solution to \eqref{3} satisfies
\begin{gather*}
\int_\mathbb R\big(u(t,x)-U(x-\gamma t)\big){\rm d}x=0,\quad t\geqslant0,\\
\lim_{t\rightarrow\infty}\|u(t,\cdot)-U(\cdot-\gamma t)\|_1=0.
\end{gather*}
\end{lemma}

\section{\Large Time-decay rate}
In this section, we will prove the estimates of time-decay rate in Theorems 1-3. The existence part in Theorem 1, however, is given by Proposition \ref{ext} and proved in the Appendix A.

\subsection{Estimates in Theorem 1}
Firstly, we need the estimate in $L^2$-norm. Let $u(t,x)$ be a weak solution to \eqref{1}, $X(t)$ with $X(0)=0$ be a Lipschitz continuous space shift to be determined below. Set
\begin{equation}\label{v0}
v(t,x)=u\big(t,x+X(t)\big),
\end{equation}
then
\begin{equation}\label{v}
\left\{\begin{aligned}
&v_t+(v-X^\prime)v_x=\big(|v_x|^{p-1}v_x\big)_x,\\
&v(0,x)=u_0(x).
\end{aligned}\right.
\end{equation}
Set $w\big(t,x\big)=v(t,x)-U(x)$, then
\begin{equation}\label{phi}
w_t+(v-X^\prime)v_x-(U-\gamma)U^\prime=\big(|v_x|^{p-1}v_x\big)_x
-\big(|U^\prime|^{p-1}U^\prime\big)^\prime
\end{equation}
Note that $v(0,\cdot)-U(\cdot)\in L^2(\mathbb R)$. Multiplying \eqref{phi} with $ w$ and integrating the resultant equation, we have
\begin{equation}\label{phi1}
\begin{aligned}
\frac{1}{2}\frac{{\rm d}}{{\rm d}t}\int_\mathbb R w(t,x)^2{\rm d}x+&\int_\mathbb R(|v_x|^{p-1}v_x-|U^\prime|^{p-1}U^\prime)(v_x-U^\prime){\rm d}x\\
&\hspace{-2cm}=\int_\mathbb R(X^\prime v_x-\gamma U^\prime)(v-U){\rm d}x+\frac{1}{2}\int_\mathbb R(v^2-U^2)(v_x-U^\prime){\rm d}x
\end{aligned}
\end{equation}
Using Lemma \ref{abm}, it holds that
\begin{equation}\label{bu1}
\int_\mathbb R(|v_x|^{p-1}v_x-|U^\prime|^{p-1}U^\prime)(v_x-U^\prime){\rm d}x\geqslant\frac{5}{6}\int_\mathbb R\big(c_0|v_x-U^\prime|^{p+1}+|U^\prime|^{p-1}(v_x-U^\prime)^2\big){\rm d}x.
\end{equation}
For the terms on the right-hand side, we have, by denoting $y=U(x)$, that
\begin{equation}\label{bu2}
\int_\mathbb R(X^\prime v_x-\gamma U^\prime)(v-U){\rm d}x=(X^\prime-\gamma)\int_{u_+}^{u_-}\!\!(v-U)U^\prime{\rm d}x=-(X^\prime-\gamma)\int_{u_+}^{u_-}\!\! w\big(t,U^{-1}(y)\big){\rm d}y
\end{equation}
and
\begin{equation}\label{bu3}
\begin{aligned}
&\int_\mathbb R(v^2-U^2)(v_x-U^\prime){\rm d}x\\
&=-\int_{u_+}^{u_-}\!\!\big(v^2-U^2-2U(v-U)\big)U^\prime{\rm d}x+\int_\mathbb R\big((v^2-U^2)v_x-2U(v-U)U^\prime\big){\rm d}x\\
&=\int_{u_+}^{u_-}\!\!\left(\Big(w\big(t,U^{-1}(y)\big)+y\Big)^2-y^2-2y w\big(t,U^{-1}(y)\big)\right){\rm d}y-2\int_\mathbb R(v^2v_x-U^2U^\prime){\rm d}x\\
&=\int_{u_+}^{u_-}\!\! w^2\big(t,U^{-1}(y)\big){\rm d}y.
\end{aligned}
\end{equation}
Thus, it follows from (\ref{phi1}-\ref{bu3}) that
\begin{equation*}
\frac{{\rm d}}{{\rm d}t}\int_\mathbb R(v-U)^2{\rm d}x+D(t)\leqslant0,
\end{equation*}
where
\begin{equation}\label{d}
\begin{aligned}
D(t)=&\frac{5}{3}\int_\mathbb R\big(c_0|v_x-U^\prime|^{p+1}+|U^\prime|^{p-1}(v_x-U^\prime)^2\big){\rm d}x\\
&\hspace{2cm}+2(X^\prime-\gamma)\int_{u_+}^{u_-}\!\! w\big(t,U^{-1}(y)\big){\rm d}y-\int_{u_+}^{u_-}\!\! w^2\big(t,U^{-1}(y)\big){\rm d}y.
\end{aligned}
\end{equation}

Now we can define the spacial shift $X$. Set $X(t)$ to be the solution to
\begin{equation}\label{x}
\left\{\begin{aligned}
&X^\prime(t)=\gamma-\frac{1}{2(u_--u_+)}\int_\mathbb R\Big(u\big(t,x+X(t)\big)-U(x)\Big)U^\prime(x){\rm d}x,\\
&X(0)=0.
\end{aligned}\right.
\end{equation}
Since
\begin{equation*}
X^\prime=\gamma-\frac{1}{2(u_--u_+)}\int_\mathbb R\Big(u(t,x)-U\big(x-X(t)\big)\Big)U^\prime\big(x-X(t)\big){\rm d}x,
\end{equation*}
we conclude the existence and uniqueness of a solution $X$ to \eqref{x} by Lipschitz continuity of $U$ and $U^\prime$ from \eqref{upr} for $1\leqslant p\leqslant2$. Moreover,
\begin{equation}\label{xprest}
\begin{aligned}
|X^\prime-\gamma|&=\frac{1}{2(u_--u_+)}\left|\int_\mathbb R(v-U)U^\prime{\rm d}x\right|\\
&\leqslant\frac{1}{2(u_--u_+)}\left(\int_\mathbb R|U^\prime|{\rm d}x\right)^\frac{1}{2}\left(\int_\mathbb R(v-U)^2|U^\prime|{\rm d}x\right)^\frac{1}{2}\\
&\leqslant(u_--u_+)^{\frac{1}{p}-\frac{1}{2}}\|v-U\|_2
\end{aligned}
\end{equation}
by using H\"{o}lder's inequality and \eqref{uprest}.

\vskip 0.2in

Set $\bar w(t)=\frac{1}{u_--u_+}\int_{u_+}^{u_-}\!\!w\big(t,U^{-1}(y)\big){\rm d}y$ to be the mean of $ w$, then
\begin{equation*}
X^\prime-\gamma=\frac{1}{2(u_--u_+)}\int_{u_+}^{u_-}\!\! w\big(t,U^{-1}(y)\big){\rm d}y=\frac{\bar w(t)}{2}
\end{equation*}
from \eqref{x}. With direct calculation, it holds
\begin{equation*}
\begin{aligned}
2(X^\prime-\gamma)\int_{u_+}^{u_-}\!\! w&\big(t,U^{-1}(y)\big){\rm d}y-\int_{u_+}^{u_-}\!\! w^2w\big(t,U^{-1}(y)\big){\rm d}y\\
&=(u_--u_+)\bar w^2-\int_{u_+}^{u_-}\!\! w^2{\rm d}y=-\int_{u_+}^{u_-}( w-\bar w)^2{\rm d}y.
\end{aligned}
\end{equation*}
Lemma 3.1 in \cite{kan17} implies
\begin{equation*}
\int_{u_+}^{u_-}( w-\bar w)^2{\rm d}y\leqslant\frac{5}{6}\int_{u_+}^{u_-}(u_--y)(y-u_+) w_y^2{\rm d}y.
\end{equation*}
On the other hand,
\begin{equation*}
\int_{u_+}^{u_-}(u_--y)(y-u_+) w_y^2{\rm d}y=2\int_\mathbb R|U^\prime|^{p-1}(v_x-U^\prime)^2{\rm d}x
\end{equation*}
from \eqref{upr}. Therefore, by \eqref{d},
\begin{equation}
D(t)\geqslant\frac{5}{3}c_0\int_\mathbb R|v_x-U^\prime|^{p+1}{\rm d}x,
\end{equation}
which means
\begin{equation}\label{est0}
\frac{{\rm d}}{{\rm d}t}\|v-U\|_2^2+\frac{5}{3}c_0\|v_x-U^\prime\|_{p+1}^{p+1}\leqslant0.
\end{equation}
Hence, $v(t,\cdot)-U\in L^2(\mathbb R)$ and
\begin{equation}\label{l2est}
\|v(t,\cdot)-U(\cdot)\|_2\leqslant\|u_0-U\|_2.
\end{equation}

\vskip 0.2in

Next, we will derive the time-decay rate of $\|v(t,\cdot)-U\|_2$. Since
\begin{equation}\label{ii}
\|v-U\|_2\leqslant C\|v_x-U^\prime\|_{p+1}^\theta\|v-U\|_1^{1-\theta}
\end{equation}
from the interpolation inequality, where $\theta=\frac{p+1}{2(2p+1)}$, we have to estimate $\|v-U\|_1$. Set $\phi(t,x)=u(t,x)-U(x-\gamma t)$, then
\begin{equation}\label{psi}
\left\{\begin{aligned}
&\phi_t+(\frac{1}{2}\phi^2+U\phi)_x
=\big(|\phi_x+U^\prime|^{p-1}(\phi_x+U^\prime)-|U^\prime|U^\prime\big)_x,\\
&\phi(0,x)=u_0-U,\qquad \lim_{x\rightarrow\pm\infty}\phi(t,x)=0.
\end{aligned}\right.
\end{equation}
Define $j_\varepsilon(z)=z(\varepsilon^2+z^2)^{-\frac{1}{2}}$ for $z\in\mathbb R,\varepsilon>0$, then $j^\prime_\varepsilon(z)=\frac{1}{\varepsilon}\big(1+(\frac{z}{\varepsilon})^2\big)^{-\frac{3}{2}}$. Set $J_\varepsilon(z)=\int_0^z\frac{\eta}{\sqrt{\varepsilon^2+\eta^2}}{\rm d}\eta$ to be the primitive of $j_\varepsilon(z)$. It is easy to see that
\begin{equation*}
\lim_{\varepsilon\rightarrow0}J_\varepsilon(z)=|z|,\quad\lim_{\varepsilon\rightarrow0}j_\varepsilon(z)={\rm sgn}\,z
\end{equation*}
uniformly in $z\in\mathbb R$. Multiplying \eqref{psi}$_1$ with $j_\varepsilon(\phi)$ implies
\begin{equation}\label{psi1}
J_\varepsilon(\phi)_t+(\frac{1}{2}\phi^2+U\phi)_x j_\varepsilon(\phi)
=\big(|\phi_x+U^\prime|^{p-1}(\phi_x+U^\prime)-|U^\prime|U^\prime\big)_x j_\varepsilon(\phi).
\end{equation}
Note that
\begin{equation*}
\begin{aligned}
(\frac{1}{2}\phi^2+\phi U)_x j_\varepsilon(\phi)&=\frac{\phi^2}{\sqrt{\varepsilon^2+\phi^2}}\phi_x
+\frac{U\phi}{\sqrt{\varepsilon^2+\phi^2}}\phi_x
+\frac{\varepsilon^2+\phi^2-\varepsilon^2}{\sqrt{\varepsilon^2+\phi^2}}U^\prime\\
&=\frac{\phi^2}{\sqrt{\varepsilon^2+\phi^2}}\phi_x+(U\sqrt{\varepsilon^2+\phi^2})_x
-\frac{\varepsilon^2}{\sqrt{\varepsilon^2+\phi^2}}U^\prime.
\end{aligned}
\end{equation*}
It follows from integrating \eqref{psi1} with respect to $\xi$ over $\mathbb R$ that
\begin{equation}\label{psi2}
\frac{{\rm d}}{{\rm d}t}\int_\mathbb RJ_\varepsilon(\phi){\rm d}x+\int_\mathbb R\big(|\phi_x+U^\prime|^{p-1}(\phi_x+U^\prime)-|U^\prime|U^\prime\big)j^\prime_\varepsilon(\phi)\phi_x{\rm d}\xi\leqslant\varepsilon(u_--u_+),
\end{equation}
since $U^\prime\leqslant0$. Thus,
\begin{equation*}
\int_\mathbb RJ_\varepsilon\big((\phi(t,\cdot)\big){\rm d}x\leqslant\int_\mathbb RJ_\varepsilon\big((\phi(0,\cdot)\big){\rm d}x+\varepsilon(u_--u_+)t
\end{equation*}
for any $t>0$. Let $\varepsilon\rightarrow0$, then we have
\begin{equation}\label{1est1}
\|u(t,\cdot)-U(\cdot-\gamma t)\|_1\leqslant\|u_0-U\|_1.
\end{equation}

We now only need the estimate of $\|U(\cdot-\gamma t)-U\big(\cdot-X(t)\big)\|_1$, which follows from the estimate of $|X(t)-\gamma t|$. The proof is same as in \cite{kan17}, but for completeness, we still state the progress here. Set $\tau=X(t)-\gamma t$ and $F(\tau)=\|U(\cdot+\tau)-U(\cdot)\|_2^2$. Obviously, $F$ is even. If $\tau>0$, then
\begin{equation*}
F^\prime(\tau)=2\int_\mathbb R\big(U(x+\tau)-U(x)\big)U^\prime(x+\tau){\rm d}x=2\int_\mathbb R\int_x^{x+\tau}\!\!U^\prime(y)U^\prime(x+\tau){\rm d}y{\rm d}x.
\end{equation*}
For $\tau\geqslant1$, it holds
\begin{equation*}
F^\prime(\tau)\geqslant2\int_\mathbb R\int_{x-1}^xU^\prime(y)U^\prime(x){\rm d}y{\rm d}x=:\lambda>0.
\end{equation*}
Hence,
\begin{equation*}
F(\tau)\geqslant F(1)+\lambda(\tau-1)\geqslant\lambda(\tau-1),
\end{equation*}
which implies
\begin{equation}\label{bu5}
|\tau|\leqslant\frac{1}{\lambda}F(\tau)+1.
\end{equation}
It is easy to see that \eqref{bu5} holds true for $0<\tau<1$. Since $F$ is even, we also have \eqref{bu5} for $\tau\leqslant0$. On the other hand,
\begin{equation*}
\begin{aligned}
\big(v(t,x)-U(x)\big)^2&=\Big(u\big(t,x+X(t)\big)-U\big(x+X(t)-\gamma t\big)+U\big(x+X(t)-\gamma t\big)-U(x)\Big)^2\\
&\geqslant2\Big(u\big(t,x+X(t)\big)-U\big(x+X(t)-\gamma t\big)\Big)\Big(U\big(x+X(t)-\gamma t\big)-U(x)\Big)\\
&\qquad\qquad+\Big(U\big(x+X(t)-\gamma t\big)-U(x)\Big)^2
\end{aligned}
\end{equation*}
Thus,
\begin{equation*}
\begin{aligned}
|X(t)-\gamma t|&\leqslant\frac{1}{\lambda}\|U\big(\cdot+X(t)-\gamma t\big)-U(\cdot)\|_2^2+1\\
&\leqslant\frac{1}{\lambda}\big(\|v-U\|_2^2+2\|U(\cdot-\gamma t)-U\big(\cdot-X(t)\big)\|_\infty\|u(t,\cdot+\gamma t)-U(\cdot)\|_1\big)+1\\
&\leqslant C(\|u_0-U\|_2^2+\|u_0-U\|_1+1),
\end{aligned}
\end{equation*}
where we used \eqref{l2est} and \eqref{1est1}. Then
\begin{equation}\label{psi12}
\begin{aligned}
&\|U(\cdot-\gamma t)-U\big(\cdot-X(t)\big)\|_1\\
&\leqslant|X(t)-\gamma t|\int_0^1\!\!\int_\mathbb R-U^\prime\big(x-\mu\gamma t-(1-\mu)X(t)\big){\rm d}x{\rm d}\mu\\
&\leqslant C(\|u_0-U\|_2^2+\|u_0-U\|_1+1).
\end{aligned}
\end{equation}
Hence, by \eqref{1est1} and \eqref{psi12}, we obtain
\begin{equation}\label{l1est}
\begin{aligned}
\|v-U\|_1&\leqslant\|u(t,\cdot)-U(\cdot-\gamma t)\|_1+\|U(\cdot-\gamma t)-U\big(\cdot-X(t)\big)\|_1\\
&\leqslant C(\|u_0-U\|_2^2+\|u_0-U\|_1+1).
\end{aligned}
\end{equation}
It follows from \eqref{est0}, \eqref{ii} and \eqref{l1est} that
\begin{equation}\label{est00}
\frac{{\rm d}}{{\rm d}t}\|v-U\|_2^2+c\big(\|v-U\|_2^2\big)^{2p+1}\leqslant0.
\end{equation}
Using Lemma \ref{yt}, we have
\begin{equation}\label{est1}
\|v-U\|_2\leqslant C(1+t)^{-\frac{1}{4p}},
\end{equation}
which is actually \eqref{re0} and hence,
\begin{equation*}
\big\|u\big(t,\cdot+X(t)\big)-U(\cdot)\big\|_{L^\infty(\mathbb R)}\leqslant C(1+t)^{-\frac{1}{2p(p+3)}}
\end{equation*}
by using the interpolation inequality and $\|u_x(t,\cdot)\|_{p+1}\leqslant C\big(\|u_0-U\|_2,\|u_{0x}\|_{p+1}\big)$ for any $t>0$ from \eqref{bu4} in Appendix A, which completes the proof of estimates in Theorem 1.

\subsection{Proof of Theorem 2}
For the case that $p=1$, we can obtain the estimate of first-order derivative. From \cite{hof85}, for every $R,\tau>0$, there exists a constant $M$ such that the solution $u$ to \eqref{lvp} with \eqref{inid} satisfies $\|u_x(t,\cdot)\|_\infty\leqslant M$ for all $t>\tau$, whenever $\|u_0\|_\infty\leqslant R$. Thus, $\lim_{x\rightarrow\pm\infty}(u-U)=0$ from $u-U\in L^2(\mathbb R)$, and hence, \eqref{bndd} holds. Therefore, \eqref{est1} holds for $p=1$. Set $\psi(t,x)=u\big(t,x+X(t)\big)-U(x)$, then
\begin{equation}\label{400}
\psi_t-X^\prime\psi_x-(X^\prime-\gamma)U^\prime+(\psi+U)(\psi_x+U^\prime)-UU^\prime=\psi_{xx}.
\end{equation}
Since $u_0\in L^1(\mathbb R)\cap L^\infty(\mathbb R)$, then $\|u(t,\cdot)\|_\infty\leqslant\|u_0\|_\infty$ for any $t>0$, and then
\begin{equation}\label{41}
\|\psi(t,\cdot)\|_\infty\leqslant\|u_0\|_\infty+\|U\|_\infty.
\end{equation}
Firstly, by a similar proof of \eqref{est0} and \eqref{est1} with $p=1$, we have
\begin{equation}\label{402}
\frac{{\rm d}}{{\rm d}t}\|\psi\|_2^2+\|\psi_x\|_2^2\leqslant0
\end{equation}
and
\begin{equation}\label{401}
\|\psi(t,\cdot)\|_2\leqslant C(1+t)^{-\frac{1}{4}}.
\end{equation}
Next, we will estimate $\|\psi\|_\infty$. Multiplying \eqref{400} by $-\psi_{xx}$ and integrating the resultant equation with respect to $x$ over $\mathbb R$, we have
\begin{equation}\label{42}
\frac{{\rm d}}{{\rm d}t}\|\psi_x\|_2^2+\frac{3}{2}\|\psi_{xx}\|_2^2\leqslant C\left(\|\psi\|_2^2+\|\psi_x\|_2^2+\|\psi_x\|_3^3\right),
\end{equation}
where we used H\"{o}lder's inequality, $\|U^\prime\|_2\leqslant C$ and \eqref{xprest}. The interpolation inequality and Young's inequality imply that
\begin{equation}\label{43}
\|\psi_x\|_3^3\leqslant\frac{1}{2}\|\psi_{xx}\|_2^2+C\|\psi\|_2^{10}.
\end{equation}
Hence,
\begin{equation}\label{44}
\frac{{\rm d}}{{\rm d}t}\big(\|\psi\|_2^2+\|\psi_x\|_2^2\big)+c\big(\|\psi_x\|_2^2+\|\psi_{xx}\|_2^2\big)\leqslant C\big(\|\psi\|_2^2+\|\psi\|_2^{10}\big)
\end{equation}
by comparing (\ref{402}-\ref{43}). On the other hand, similar to the proof of \eqref{1est1}, we can obtain that $\|\psi(t,\cdot)\|_1\leqslant\|u_0-U\|_1$. Then the interpolation inequality implies
\begin{equation*}
\|\psi\|_2\leqslant C\|\psi\|_1^\frac{2}{3}\|\psi_x\|_2^\frac{1}{3}\leqslant C\|\psi_x\|_2^\frac{1}{3},\quad \|\psi_x\|_2\leqslant C\|\psi\|_\infty^\frac{2}{3}\|\psi_{xx}\|_2^\frac{1}{3}\leqslant C\|\psi_{xx}\|_2^\frac{1}{3}
\end{equation*}
by using \eqref{41}. Thus, from \eqref{401} and \eqref{44}, it holds
\begin{equation*}
\frac{{\rm d}}{{\rm d}t}\big(\|\psi\|_2^2+\|\psi_x\|_2^2\big)+c\big(\|\psi\|_2^2+\|\psi_x\|_2^2\big)^3\leqslant C(1+t)^{-\frac{1}{2}}.
\end{equation*}
Using Lemma \ref{yt}, we have
\begin{equation*}
\|\psi\|_2^2+\|\psi_x\|_2^2\leqslant C(1+t)^{-\frac{1}{6}}
\end{equation*}
for suitably large $t>\tau$, which, together with \eqref{41}, implies
\begin{equation}\label{45}
\|\psi(t,\cdot)\|_\infty\leqslant C\|\psi\|_2^\frac{1}{2}\|\psi_x\|_2^\frac{1}{2}\leqslant C(1+t)^{-\frac{1}{6}}.
\end{equation}
The proof is completed.
\begin{rem}
In fact, the discussion above can be extended to the solution to Cauchy problem
\begin{equation}\label{bu6}
\left\{\begin{aligned}
&u_t+f(u)_x=u_{xx},\\
&u(0,\xi)=u_0(x)
\end{aligned}\right.
\end{equation}
with
\begin{equation*}
\sup_{u\in\mathbb R}\left|f(u)-\frac{1}{2}u^2\right|<\frac{1}{11},
\end{equation*}
since the $L^2$ estimate has be obtained in \cite{kan17}. Thus, the estimates in Theorem 2 hold true for \eqref{bu6} with $u_0\in L^1(\mathbb R)\cap L^\infty(\mathbb R)$.
\end{rem}

\subsection{Brief proof of Theorem 3}
For the case with general convex flux, since $u_0-U\in L^1(\mathbb R)$, there exists a space shift $y$ such that $\int_\mathbb R\big(u_0(x)-U(x-y)\big){\rm d}x=0$. Without loss of generality, we suppose $y=0$ in the rest of this section.

Set
\begin{gather*}
\phi(t,\xi)=u(t,\xi+\gamma t)-U(\xi),\quad \Phi(t,\xi)=\int_{-\infty}^\xi\phi(t,\eta){\rm d}\eta,\\
\phi_0(\xi)=u_0(\xi)-U(\xi),\quad \Phi_0(\xi)=\int_{-\infty}^\xi\phi_0(\eta){\rm d}\eta.
\end{gather*}
From \cite{hof85} we obtain the existence, uniqueness and regularity of the solution to \eqref{3}. In addition, for every $R,\tau>0$, there exists a constant $M$ such that the solution $u$ to \eqref{lvp} with \eqref{inid} satisfies $\|u_x(t,\cdot)\|_\infty\leqslant M$ for all $t>\tau$, whenever $\|u_0\|_\infty\leqslant R$. Thus, $\lim_{x\rightarrow\pm\infty}(u-U)=0$ from $u-U\in L^2(\mathbb R)$. We then have
\begin{equation}\label{50}
\left\{\begin{aligned}
&\Phi_t-\gamma\Phi_\xi+f(U+\Phi_\xi)-f(U)=\Phi_{\xi\xi},\\
&\Phi(0,\xi)=\Phi_0(\xi),\\
&\lim_{\xi\rightarrow\pm\infty}\Phi(t,\xi)=0
\end{aligned}\right.
\end{equation}
and
\begin{equation}\label{500}
\left\{\begin{aligned}
&\phi_t-\gamma\phi_\xi+\big(f(U+\phi)-f(U)\big)_\xi=\phi_{\xi\xi},\\
&\phi(0,\xi)=\phi_0(\xi),\\
&\lim_{\xi\rightarrow\pm\infty}\phi(t,\xi)=0.
\end{aligned}\right.
\end{equation}
Multiplying \eqref{50}$_1$ by $|\Phi|^{r-2}\Phi$ with $r\geqslant2$, we have
\begin{equation}\label{51}
\begin{aligned}
\frac{{\rm d}}{{\rm d}t}\left(\frac{1}{r}|\Phi|^r\right)+&(r-1)|\Phi|^{r-2}\Phi_\xi^2-\frac{1}{r}f^{\prime\prime}(U)U^\prime|\Phi|^r\\
&=(\cdots)_\xi-\frac{1}{2}f^{\prime\prime}(U+\theta\Phi_\xi)|\Phi|^{r-2}\Phi\Phi_\xi^2,
\end{aligned}
\end{equation}
where $\theta\in[0,1]$ and
\begin{equation*}
(\cdots)=\frac{\gamma}{r}|\Phi|^r-\frac{1}{r}f^\prime(U)|\Phi|^r+|\Phi|^{r-2}\Phi\Phi_\xi.
\end{equation*}
Noting that $U^\prime<0,f^{\prime\prime}>0$ and $\|\Phi\|_\infty\leqslant\|\phi\|_1\rightarrow0$ as $t\rightarrow\infty$ by Lemma \ref{phil1est}, it follows from \eqref{51} for $t$ sufficiently large that
\begin{equation}\label{52}
\frac{{\rm d}}{{\rm d}t}\|\Phi\|_r^r+c\int_\mathbb R|\Phi|^{r-2}\Phi_\xi^2{\rm d}\xi\leqslant0.
\end{equation}
From \eqref{52} we can obtain the decay rate of $\|\Phi\|_r$. Using this rate in energy estimates of $\phi$ and $\phi_\xi$ implies \eqref{bu10}, since $\phi(t,\cdot)\in L^\infty(\mathbb R)$ for any $t\geqslant0$ and $\phi_\xi(t,\cdot)\in L^\infty(\mathbb R)$ for suitably large $t$. The details are the same as in \cite{hua22}, so we omit it here.

\section{\Large Appendix}
\subsection*{A: Proof of Proposition \ref{ext}}
Since $u_0-U\in L^2(\mathbb R)$ and $u_{0x}\in L^{p+1}(\mathbb R)$, we can construct a sequence of smooth functions $\{u_0^\varepsilon\}$, which uniformly converges to $u_0$ as $\varepsilon\rightarrow0+$, and satisfies
\begin{equation*}
\|u_0^\varepsilon-U\|_2\leqslant\|u_0-U\|_2,\quad \|u_{0\xi}^\varepsilon\|_{p+1}\leqslant(1+K\varepsilon)\|u_{0\xi}\|_{p+1},
\end{equation*}
for some constant $K>0$ depending only on the approximation process. Consider the Cauchy problem
\begin{equation}\label{a0}
\left\{\begin{aligned}
&u_t^\varepsilon-\gamma u^\varepsilon_\xi+u^\varepsilon u^\varepsilon_\xi=\left(\big((u_\xi^\varepsilon)^2+\varepsilon\big)^\frac{p-1}{2}u_\xi^\varepsilon\right)_\xi,\\
&u^\varepsilon(0,\xi)=u_0^\varepsilon(\xi).
\end{aligned}\right.
\end{equation}

\vskip 0.2in

Let $\phi(t,\xi)=u^\varepsilon(t,\xi)-U(\xi)$, where $U(\xi)$ is given by \eqref{U}. It is easy to see that
\begin{equation}\label{a1}
\phi_t-\gamma\phi_\xi+(\phi+U)(\phi_\xi+U^\prime)-UU^\prime
=\left(\big((\phi_\xi+U^\prime)^2+\varepsilon\big)^\frac{p-1}{2}(\phi_\xi+U^\prime)-|U^\prime|^{p-1}U^\prime\right)_\xi.
\end{equation}
Multiplying \eqref{a1} by $\phi$ implies
\begin{equation}\label{a2}
\begin{aligned}
\frac{1}{2}(\phi^2)_t+&(\cdots)_\xi
+\big(|\phi_\xi+U^\prime|^{p-1}(\phi_\xi+U^\prime)-|U^\prime|^{p-1}U^\prime\big)\phi_\xi\\
&=(\phi^2+2\phi U)\phi_\xi-\left(\big((\phi_\xi+U^\prime)^2+\varepsilon\big)^\frac{p-1}{2}
-|\phi_\xi+U^\prime|^{p-1}\right)(\phi_\xi+U^\prime)\phi_\xi,
\end{aligned}
\end{equation}
where
\begin{equation*}
(\cdots)=-\frac{\gamma}{2}\phi^2+\frac{1}{2}\big((\phi+U)^2-U^2\big)\phi
-\left(\big((\phi_\xi+U^\prime)^2+\varepsilon\big)^\frac{p-1}{2}(\phi_\xi+U^\prime)-|U^\prime|^{p-1}U^\prime\right)\phi.
\end{equation*}
It is noted that for any $a,b\in\mathbb R,p\geqslant1$, the inequality
\begin{equation}\label{ab1}
\big(|a|^{p-1}a-|b|^{p-1}b\big)(a-b)\geqslant c_p|a-b|^{p+1}
\end{equation}
holds true, and for any $a,b\geqslant0,0\leqslant q\leqslant1$,
\begin{equation}\label{ab2}
|a^q-b^q|\leqslant c_p|a-b|^q.
\end{equation}
In addition, if $p>3$,
\begin{equation}\label{ab3}
\begin{aligned}
\big((\phi_\xi+U^\prime)^2+\varepsilon\big)^\frac{p-1}{2}-|\phi_\xi+U^\prime|^{p-1}
&=\frac{p-1}{2}\varepsilon\big(\theta\varepsilon+(\phi_\xi+U^\prime)^2\big)^\frac{p-3}{2}\\
&\leqslant c|\phi_\xi|^{p+1}+C\varepsilon^\kappa
\end{aligned}
\end{equation}
for some $\kappa>0$ by using Young's inequality and \eqref{uprest}. Then, integrating \eqref{a2} with respect to $\xi$ over $\mathbb R$ and using (\ref{ab1}-\ref{ab3}), we have
\begin{equation}\label{a3}
\frac{{\rm d}}{{\rm d}t}\|\phi\|_2^2+c\|\phi_\xi\|_{p+1}^{p+1}\leqslant C\big(\|\phi\|_2^2+\varepsilon^\nu\big)
\end{equation}
for some $\nu>0$ by choosing $\varepsilon$ suitably small, where we used Young's inequality on the right-hand side. Thus, for any given $T>0$, it follows
\begin{equation}\label{a4}
\|\phi(T,\cdot)\|_2^2+c\int_0^T\|\phi_\xi(t,\cdot)\|_{p+1}^{p+1}{\rm d}t\leqslant e^{CT}\|\phi(0,\cdot)\|_2^2+C\varepsilon^\nu T.
\end{equation}
Next, similar to the progress in \cite{mat94}, it holds
\begin{gather*}
\|u_\xi^\varepsilon(T,\cdot)\|_{p+1}^{p+1}+\varepsilon^\frac{p-1}{2}\|u_\xi^\varepsilon(T,\cdot)\|_2^2
+\int_0^T\!\!\int_\mathbb R\big((u_\xi^\varepsilon)^2+\varepsilon\big)^{p-1}(u_{\xi\xi}^\varepsilon)^2{\rm d}\xi{\rm d}t\leqslant C,\\
\int_0^T\!\!\int_\mathbb R\big((u_\xi^\varepsilon)^2+\varepsilon\big)^\frac{p-1}{2}|u_\xi^\varepsilon|^3{\rm d}\xi{\rm d}t\leqslant C.
\end{gather*}
Here, $C=C\big(\|u_0^\varepsilon-U\|_2,\|u_{0\xi}^\varepsilon\|_{p+1},\varepsilon\|u_{0\xi}^\varepsilon\|_2,\varepsilon T\big)$. Then, let $\varepsilon\rightarrow0$, it holds
\begin{equation}\label{bu4}
\|u_\xi(T,\cdot)\|_{p+1}^{p+1}\leqslant C\big(\|u_0-U\|_2,\|u_{0\xi}\|_{p+1}\big)
\end{equation}
for any $T>0$. Since $u_0-U\in L^2(\mathbb R)$ and $u_{0\xi}\in L^{p+1}(\mathbb R)$, we obtain that \eqref{1} has a weak solution $u(t,\xi)$ on $[0,T]$ for any $T>0$, and
\begin{gather*}
u-U\in C\big([0,\infty);L^2(\mathbb R)\big)\cap L^\infty\big([0,\infty);L^2(\mathbb R)\big),\\
u_\xi\in L^\infty\big([0,\infty);L^{p+1}(\mathbb R)\big)\cap L^{p+1}\big([0,\infty)\times\mathbb R\big),\\
\big(|u_\xi|^{p-1}u_\xi\big)_\xi\in L^2\big([0,\infty)\times\mathbb R\big).
\end{gather*}

\vskip 0.2in

It remains to check the uniqueness. Suppose $u,v$ are weak solutions to \eqref{1} with initial data $u_0=v_0$ and let $\psi=u-v$. Similar to the discussion about \eqref{a4} above, we have
\begin{equation}
\|\psi(T,\cdot)\|_2^2+c\int_0^T\|\psi_\xi(t,\cdot)\|_{p+1}^{p+1}{\rm d}t\leqslant e^{CT}\|\psi(0,\cdot)\|_2^2=0
\end{equation}
for any $T>0$. Then, $u=v$ and the proof is completed.

\vskip 0.3in

\subsection*{B: Proof of Lemma \ref{abm}}
Obviously, \eqref{ab} holds true if $a=b$ or $p=1$ by choosing $c_0\leqslant\frac{2^{1-p}}{5}$. Without loss of generality, we suppose $a>b$ and $p>1$ from now on. Set
\begin{equation*}
I=\frac{(|a|^{p-1}a-|b|^{p-1}b)(a-b)}{\big(c_0|a-b|^{p-1}+(\max\{|a|,|b|\})^{p-1}\big)(a-b)^2}.
\end{equation*}

\vskip 0.2in

\noindent{\bf Case 1} $(ab\geqslant0)$. Set
\begin{equation*}
\theta=\left\{\begin{aligned}
&\frac{b}{a},\quad 0\leqslant b<a,\\
&\frac{a}{b},\quad b<a\leqslant0,
\end{aligned}\right.
\end{equation*}
then $0\leqslant\theta<1$ and
\begin{equation*}
I=\frac{1-\theta^p}{\big(c_0(1-\theta)^{p-1}+1\big)(1-\theta)}=:h_1(\theta).
\end{equation*}
It is easy to see that
\begin{equation*}
h_1(0)=\frac{1}{c_0+1}\geqslant\frac{5}{6},\quad \lim_{\theta\rightarrow1-}h_1(\theta)=p>\frac{5}{6}
\end{equation*}
by choosing $c_0\leqslant\frac{1}{5}$. On the other hand,
\begin{equation*}
h_1^\prime(\theta)=\frac{c_0p(1-\theta^{p-1})(1-\theta)^{p-1}
+1-p\theta^{p-1}+(p-1)\theta^p}{\big(c_0(1-\theta)^{p-1}+1\big)^2(1-\theta)^2}>0.
\end{equation*}
Thus, $h_1(\theta)\geqslant\frac{5}{6}$ on $[0,1)$. Hence, \eqref{ab} holds true.

\vskip 0.2in

\noindent{\bf Case 2} $(a>0>b)$. By direct calculation, we have
\begin{equation*}
I=\frac{|a|^p+|b|^p}{\big(c_0(|a|+|b|)^{p-1}+(\max\{|a|,|b|\})^{p-1}\big)(|a|+|b|)}.
\end{equation*}
Without loss of generality, we soppose $a\geqslant-b$. Set $\theta=-\frac{b}{a}$, then $0<\theta\leqslant1$ and
\begin{equation*}
I=\frac{1+\theta^p}{\big(c_0(1+\theta)^{p-1}+1\big)(1+\theta)}=:h_2(\theta).
\end{equation*}
It is easy to see that
\begin{equation*}
\lim_{\theta\rightarrow0+}h_2(\theta)=\frac{1}{c_0+1}>\frac{5}{6},\quad h_2(1)=\frac{1}{2^{p-1}c_0+1}\geqslant\frac{5}{6}
\end{equation*}
by choosing $c_0\leqslant\frac{2^{1-p}}{5}$. Since $h_2$ is smooth, we only need to consider the extreme value of $h_2$ on $(0,1)$. Let $h_2^\prime=0$, we have
\begin{equation*}
0=c_0p(1-\theta^{p-1})(1+\theta)^{p-1}+1-p\theta^{p-1}-(p-1)\theta^p=:g(\theta).
\end{equation*}
Obviously, $\lim_{\theta\rightarrow0+}g(\theta)=c_0p+1>0$, $\lim_{\theta\rightarrow1-}g(\theta)=-2(p-1)<0$ and
\begin{equation*}
g^\prime(\theta)=p(p-1)\big(c_0(1-\theta^{p-2}-2\theta^{p-1})(1+\theta)^{p-2}-\theta^{p-2}-\theta^{p-1}\big).
\end{equation*}
Since $1<p<2$, it holds $g^\prime<0$ by choosing $c_0<1$. Hence, there uniquely exists a $\theta_p\in(0,1)$ such that $g(\theta_p)=0$ and $h_2(\theta_p)=\min_{0<\theta\leqslant1}h_2(\theta)$.

To ensure $h_2\geqslant\frac{5}{6}$, we need to consider the extreme case. Denote $\tilde h_2(\theta)=\frac{1+\theta^p}{1+\theta}$ to be $h_2$ with $c_0=0$. Since for any given $\theta\in(0,1)$, $\tilde h_2(\theta)$ is monotone decreasing with respect to $m$. Then, there exists a constant $p_0\in(\frac{39}{20},\frac{59}{30})$ such that, for any $p\in[1,p_0)$ and $\theta\in[0,1]$, $\tilde h_2>\frac{5}{6}+\delta$ with $\delta>0$ being a small constant. Hence, there exists a constant $c_0>0$ such that $h_2\geqslant\frac{5}{6}$ for any $m\in[1,p_0)$ and $\theta\in[0,1]$. Thus, the proof is completed.

\end{document}